\documentclass{birkjour}

\usepackage{subfigure}
 \newtheorem{theorem}{Theorem}[section]
 \newtheorem{corollary}[theorem]{Corollary}
 \newtheorem{lemma}[theorem]{Lemma}
 \newtheorem{proposition}[theorem]{Proposition}

 \numberwithin{equation}{section}
 
 \newcommand{\R}{\mathbb{R}}

  \newcommand{\e}{\mathbf{e}}
  \newcommand{\h}{\mathbf{h}}
   
    \newcommand{\n}{\mathbf{n}}

  \renewcommand{\v}{\mathbf{v}}

 \newcommand{\fsep}{\hspace*{\fill}}

\begin{document}

%
%
%
%
%
%
%
%
%

\title[Convex Equal-Area Polygons]
  {\centerline{Affine Properties of Convex Equal-Area Polygons}}

\author[M.Craizer]{Marcos Craizer}

\address{%
Departamento de Matem\'{a}tica- PUC-Rio\br
Rio de Janeiro\br
Brazil}

\email{craizer@puc-rio.br}

\thanks{The first and second authors want to thank CNPq for financial support during the preparation of this manuscript.}
\author[R.C.Teixeira]{Ralph C.Teixeira}
\address{Departamento de Matem\'{a}tica Aplicada- UFF \br
Niter\'oi\br
Brazil}
\email{ralph@mat.uff.br}

\author[M.A.H.B.daSilva]{Moacyr A.H.B.da Silva}
\address{Centro de Matem\'{a}tica Aplicada- FGV \br
Rio de Janeiro\br
Brazil}
\email{moacyr@fgv.br}


\date{January 24, 2010}

\begin{abstract}
In this paper we discuss some affine properties of convex equal-area polygons, which are convex polygons such that all triangles formed by three consecutive vertices have the same area.  
Besides being able to approximate closed convex smooth curves almost uniformly with respect to affine length, convex equal-area polygons admit natural definitions of the usual 
affine differential geometry concepts, like affine normal and affine curvature. These definitions lead to discrete analogous of the six vertices theorem and an affine isoperimetric inequality.
One can also define discrete counterparts of the affine evolute, parallels and the affine distance symmetry set preserving many of the properties valid for smooth curves. 
\end{abstract}

\maketitle


\section{Introduction}

A convex equal-area polygon is a simple polygon bounding a convex planar domain such that all triangles formed by three consecutive vertices have the same area (\cite{Harel03}). As a consequence, each four consecutive vertices form a trapezium.
It is worthwhile to mention that the class of equal-area polygons includes also non-convex polygons 
and appears as a set of maximizers of some area functionals (see \cite{Gronchi08}).  
{\it Affinely regular polygons} are affine transforms of convex regular polygons and thus are convex equal-area polygons. 
The affinely regular polygons appear as maximizers of area functionals and in many other different contexts (see \cite{Fischer98}).

We begin with the study of the space ${\mathcal T}_n$ of convex equal area polygons with $n$ vertices, which can be thought as a subset of the $(2n)$-dimensional affine space.
The space ${\mathcal T}_n$ has dimension $(n-5)$ and in section \ref{section:TrapezoidalPolygons} we describe in detail its structure. 
We also consider the problem of approximating a closed convex curve  by polygon in ${\mathcal T}_n$. We show how to
construct a convex equal-area $n$-gon whose vertices, for $n$ large, are close to a uniform distribution of points on the curve with respect to affine arc-length.

For convex equal-area polygons, we define in section \ref{section:DiscreteAffineGeometry} the affine tangent vector, affine normal vector at each vertex and the affine curvature $\mu$ at each edge in a straightforward way. 
In this context, we define a sextactic edge as an edge where $\Delta\mu$ is changing sign. Then we show that there are at least six sextactic edges in a convex equal-area polygon, 
which is a discrete analog of the well-known six vertices theorem.

In section \ref{section:Evolutas} we define the affine evolute and the $\lambda$-parallel curves for convex equal-area polygons. As in the smooth case, the six vertices theorem may be re-phrased as the existence of at least six cusps on the affine evolute (\cite{Ovsienko01},\cite{Tabach00}). Besides, we show that cusps of the parallels belong to the evolute, as in the smooth case.
The set of self-intersections of the parallels is called the affine distance symmetry set (ADSS) (\cite{Giblin98}). Again we prove here two properties analogous to the smooth case: cusps of the ADSS occur at points of the
affine evolute and endpoints of the ADSS are exactly the cusps of the evolute (\cite{Giblin08}).

Finally in section \ref{section:Isoperimetric} we prove an isoperimetric inequality by using a Minkowski mixed area inequality for parallel polygons. This inequality has a well-known smooth counterpart 
 saying that the integral of the affine curvature with respect to affine arc-length is at most $\frac{L^2}{2A}$, where $L$ is the affine perimeter
of the convex curve and $A$ is the area that it bounds, with equality holding only for ellipses (see \cite{Sapiro94}). In the convex equal-area polygon case, the equality holds only for affinely regular polygons.

 We have used the free software GeoGebra (\cite{GeoGebra}) for all figures and many experiments during the preparation of the paper. Applets of some of these experiments are available at \cite{Applets}. 
 We would like to thank 
 the GeoGebra team for this excellent mathematical tool. 
 
 \section{Review of some concepts of affine geometry of convex curves}
 
 In this section we review the basic affine differential concepts associated with convex curves. Since in this paper we deal 
 with a discrete model, these concepts are not strictly necessary for understanding the paper. Nevertheless, 
 they are the inspiration for most definitions and results. 
 
 \paragraph*{Affine arc-length, tangent, normal and curvature} Let $\gamma:[a,b]\to\R^2$ be a regular parameterized strictly convex curve, i.e., 
 $\left[\gamma'(t), \gamma''(t) \right]\neq 0$. The parameter $s$ defined by 
 $s(t)=\int_a^t\left[\gamma'(t), \gamma''(t) \right]^{1/3}dt\ $
is called the affine arc-length parameter. If $\gamma$ is parameterized by affine arc-length, 
\begin{equation}\label{affinearclength}
\left[\gamma'(s), \gamma''(s) \right]=1.
\end{equation}
The vector $\gamma'(s)$ is called the affine tangent and $\gamma''(s)$ is called the affine normal.
Differentiating \eqref{affinearclength} we obtain that $\gamma'''(s)$ is parallel to $\gamma'(s)$.
Thus we write
 $\gamma'''(s)+\mu(s)\gamma'(s)=0$,
 where the scalar function $\mu(s)$ is called the affine curvature. A {\it vertex} or {\it sextactic point} of $\gamma$ is a point 
 where the derivative of the curvature is zero. 
 
 \paragraph*{Affine evolute, parallels and affine distance symmetry set}
 The envelope of the family of affine normal lines is called the affine evolute of $\gamma$. For a fixed $\lambda\in\R$, the set
of points of the form 
$\gamma(s)+\lambda\gamma''(s)$, $s\in[a,b]$, 
 is called the $\lambda$-parallel. The ADSS is the locus of self-intersection of parallels (\cite{Giblin98},\cite{Giblin08}). 
 
 \paragraph*{Six vertices theorem}
 A well-known theorem of global affine planar geometry says that any convex closed planar curve must have at least six vertices (\cite{Buchin83}). 
 Equivalently, the affine evolute must have at least six cusps. 
 
 \paragraph*{An affine isoperimetric inequality} There are some different inequalities referred as affine isoperimetric inequalities. In this paper
 we shall consider the following one:
 for a closed convex curve $\gamma$ with affine perimeter $L$ enclosing a region of area $A$, 
 \begin{equation}\label{SmoothIsoIneq}
 \int_{\gamma}\mu ds\leq \frac{L^2}{2A},
\end{equation}
with equality if and only if $\gamma$ is an ellipsis  (\cite{Sapiro94}).

\section{Convex Equal-Area Polygons}\label{section:TrapezoidalPolygons}

\subsection{Basic definitions}

Consider a closed convex polygon $P$ with vertices $P_1,P_2,...,P_n$, $n\geq 5$, and oriented sides 
$$
\v_{i+\frac{1}{2}}=P_{i+1}-P_{i}, 
$$
$1\leq i\leq n$ (throughout the paper, all indices are taken modulo $n$). The polygon $P$ is said to be {\it equal-area} if 
\begin{equation}\label{DefEqualArea}
[\v_{i-\frac{1}{2}}, \v_{i+\frac{1}{2}}]=[\v_{i+\frac{1}{2}}, \v_{i+\frac{3}{2}}],
\end{equation}
for every $1\leq i\leq n$. 
This is equivalent to say that  $P_{i+2}-P_{i-1}$ is parallel to $\v_{i+\frac{1}{2}}$, for every $1\leq i\leq n$ (see figure \ref{sides10}).
We write
\begin{equation}\label{RatioDiagonalSide}
P_{i+2}-P_{i-1}=(3-\mu_{i+\frac{1}{2}})\v_{i+\frac{1}{2}},
\end{equation}
for some scalar $\mu_{i+\frac{1}{2}}<2$.

The following lemma will be useful below:

\begin{lemma}\label{Hyperbola}
Consider a trapezium $ABCD$ with $AD\parallel BC$. Then $A$ and $D$ belong to a hyperbola whose asymptotes are the lines parallel to $AB$ and $CD$ passing through $C$ and $B$, respectively.
As a consequence, in a convex equal-area polygon, $P_{i}$ and $P_{i+3}$ belong to a hyperbola whose asymptotes are $P_ {i-1}P_{i+2}$ and $P_{i+1}P_{i+4}$, for any $1\leq i\leq n$.
\end{lemma}

\begin{figure}[htb]
\centering
\includegraphics[width=0.45\linewidth]{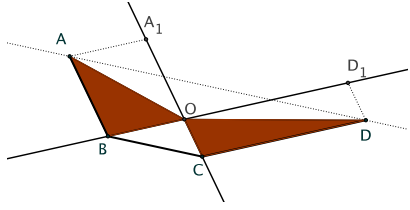}
\caption{$A$ and $D$ belong to a hyperbola with asymptotes $BO$ and $CO$.}
\label{TwoPoints}
\end{figure}

\begin{proof}
Denote by $O$ be the intersection of the proposed asymptotes,  by $A_1$ the intersection of the parallel to $CD$ through $A$ with $CO$ and by $D_1$ the intersection of the parallel 
to $AB$ through $D$ with $BO$ (see figure \ref{TwoPoints}). We must prove that the area of the parallelograms $AA_1OB$ and $DD_1OC$ area equal. But 
$$
Area(AA_1OB)=2Area(ABO), \ \ \ Area(DD_1OC)=2Area(DOC)
$$
and 
$$
Area(ABO)=Area(ABC)=Area(DBC)=Area(DOC),
$$
thus proving the lemma. 
\end{proof}

\subsection{ The space of convex equal-area polygons}

Denote by ${\mathcal T}_n$ the space of convex equal-area $n$-gons modulo affine equivalence. By taking the coordinates of the vertices of the polygon, 
we can consider ${\mathcal T}_n\subset\R^{2n}$. For a polygon to be equal area, conditions \eqref{DefEqualArea} must be satisfied. Since these
conditions are cyclic, one of them is redundant, and thus they amount to $(n-1)$ equations. Also, the affine group is $6$-dimensional and thus we expect ${\mathcal T}_n$ to be a $(n-5)$-dimensional space. In fact, this  is proved in  \cite{Harel03} for the space of (not necessarily convex) equal-area polygons. 

In this section we study the space  ${\mathcal T}_n$ in some detail. 
We prove here the following proposition:

\begin{proposition}\label{CEAspace}
There exists an open set $U\subset\R^{n-5}$ and $W\subset\partial U$ with injective smooth maps $\phi_1:U\cup W\to{\mathcal T}_n$, $\phi_2:U\cup W\to{\mathcal T}_n$ such that 
$Im(\phi_1)\cup Im(\phi_2)={\mathcal T}_n$ and $\phi_1|_W=\phi_2|_W$. 
Moreover there exists an open set $U_0\subset U$, $\overline{U_0}\cap W=\emptyset$, 
such that $\left.\phi_1\right|_{U_0}=\left.\phi_2\right|_{U_0}$.
\end{proposition}

\begin{proof}
Begin with three non-collinear points $P_1,P_2, P_3$.  Since we are working modulo affine equivalence, we may assume they are fixed. Let $\mu=(\mu_{2+\frac{1}{2}},...,\mu_{n-4+\frac{1}{2}})\in\R^{n-5}$.
Since $P_4$ must be in a line parallel to $P_2P_3$ passing through $P_1$, so it is determined by the choice of $\mu_{2+\frac{1}{2}}$. In other words, we use \eqref{RatioDiagonalSide} with $i=2$ to define $P_4$. 
Similarly, we determine 
$P_5,P_6,....,P_{n-2}$. 
Now that we know the first $(n-2)$ vertices $\{P_1,P_2,....,P_{n-2}\}$, we want to find the vertices $P_{n-1}$ and $P_n$ to complete the convex equal-area $n$-gon.
Denote by $r_1$ the line through $P_{n-2}$ and $P_{1}$, by $r_2$ the line parallel to $P_{n-3}P_{n-2}$ passing through $P_{n-4}$ and by $r_3$ the line parallel to $P_1P_2$ passing through $P_3$ 
(see figure \ref{sides10} and applet \cite{Applets}). 
The following five conditions must be satisfied:
$$
P_{n-1}P_n\parallel r_1,\ \ P_{n-1}\in r_2,\ \ P_n\in r_3,\ \ P_{n-2}P_{n-1}\parallel P_{n-3}P_{n},\ \ P_{n-1}P_{2}\parallel P_nP_1.
$$
The fifth condition is redundant, since the $n$ equations \eqref{DefEqualArea} are cyclic. Thus we have to verify only the first four.
According to lemma \ref{Hyperbola}, the first and fourth conditions say that $P_n$ belongs to the hyperbola passing through $P_{n-3}$ with asymptotes $r_1$ and $r_2$.
Thus the points $P_n$ and $P'_n$ that close the convex equal-area polygon
are the intersections of the branch $h$ of this hyperbola that does not contain $P_{n-3}$ with the line $r_3$, when these intersections exist. 

\begin{figure}[htb]
 \centering
 \includegraphics[width=0.6\linewidth]{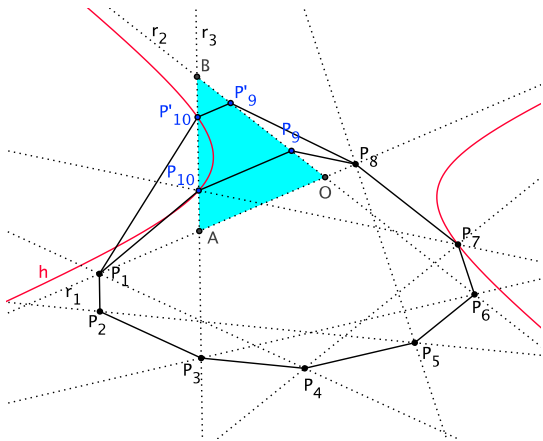}
 \caption{Conditions for the existence of a convex equal-area $n$-gon. }
\label{sides10}
\end{figure}

Let $O=r_1\cap r_2$, $A=r_1\cap r_3$ and $B=r_2\cap r_3$. 
It is clear that $A$ must be on the segment
$P_{1}O$.
We consider some different cases: (0) $B$ is in the half-line $AP_3$. In this case the branch of hyperbola $h$ intersects $r_3$ exactly once. 
(1) $r_2$ is parallel to $r_3$. In this case the branch of hyperbola $h$ also intersects $r_3$ exactly once. 
(2) When $B$ is on the half-line $P_3A$, the branch of hyperbola $h$ touches $r_3$ if and only if the ratio 
\begin{equation}\label{AreaCondition}
R(\mu)=\frac{Area(AOB)}{Area(P_1P_2P_3)} 
\end{equation}
is greater than or equal to $4$. This is not difficult to verify by arguments similar to the ones used in lemma \ref{Hyperbola}. (2a) When $R(\mu)$ is strictly greater than $4$, 
the segment $AB$ intersect $h$ is two points (see figure \ref{sides10}). 
(2b) When this ratio is exactly $4$, the segment $AB$ is tangent to $h$ and thus intersect it in one point.
We shall denote by $U_0$, $U_1$, $U_2$ and $W$ the subsets of $\mu\in\R^{n-5}$ such that (0), (1), (2a) or (2b)  occurs, respectively.

It is now clear how to define the smooth maps $\phi_1$ and $\phi_2$. For $\mu\in U_0\cup U_1\cup W$, $\phi_1=\phi_2$ correspond to the unique convex equal area polygon
corresponding to $\mu$.  For $\mu\in U_2$, let $\phi_1$ and $\phi_2$ correspond to the two convex equal-area polygons associated with $\mu$. To conclude the proof of the proposition, 
take $U=U_0\cup U_1\cup U_2$. 
\end{proof}

Next lemma says that the affinely regular $n$-gon corresponds to $\mu\in U$:

\begin{lemma}\label{AreasNgon}
Consider an affine regular $n$-gon: 
\begin{itemize}
\item For $n=5,6,7,8$, $\mu\in \overline{U_0}$.
\item For $n\geq 9$, 
the ratio $R(\mu)$ defined by equation \eqref{AreaCondition} 
is strictly greater than $4$ and thus $\mu\in U-\overline{U_0}$.
\end{itemize}
\end{lemma}

\begin{proof}
The first item is immediate. For the second item, one can obtain 
$$
R=\frac{(X+Y)^2}{XY},
$$
where $X=\frac{1}{2\cos(\alpha_n)}$, $Y=\frac{1}{\sqrt{2(1+\cos(\alpha_n))}}$ and $\alpha_n=\frac{4\pi}{n}$. 
We conclude that $R\geq 4$, and since $\cos{\alpha_n}<1$, the inequality is strict.
\end{proof}

We conjecture that the set ${\mathcal T}_n$ is connected, as suggested by  experiments done with GeoGebra. 

\subsection{ Approximating a convex curve by a convex equal-area polygon}

In this section we propose an algorithm for approximating a closed convex curve by convex equal-area polygons. Although the resulting polygons may be neither inscribed nor circumscribed, 
they are asymptotically close to the inscribed polygon whose vertices are equally spaced with respect to the affine arc-length of the curve.  For asymptotic results concerning optimal approximations
of convex curves by inscribed or circumscribed polygons, we refer to \cite{Ludwig94} and \cite{McClure75}. For surveys on approximation of convex curves by polygons,  
we refer to \cite{Gruber83} and  \cite{Gruber93}.

Given a closed convex planar curve $\gamma$, we consider the following algorithm: Fix any three points $P_1,P_2,P_3$ in a positive orientation at $\gamma$
such that $d(P_1,P_2)=d(P_2,P_3)=\frac{L}{n}$, where $d(\cdot,\cdot)$ denotes affine distance along $\gamma$ and $L$ is the affine perimeter of $\gamma$. Then $P_4$ is obtained as the intersection of a parallel to $P_2P_3$ at $P_1$ with $\gamma$. Proceeding in this way we obtain $P_k$, $k=1,...,m$ at the curve, and we continue until condition (2) in the proof of proposition \ref{CEAspace} holds with $R(\mu)<4$. 

Denote by $o(s^k)$ any quantity satisfying 
$\lim_{s\to 0}\frac{o(s^k)}{s^{k}}=0$.

\begin{lemma}
Assume that $d(P_{i-1},P_{i})=s$ and $d(P_{i},P_{i+1})=s+o(s^2)$. Then $d(P_{i+1},P_{i+2})=s+o(s^2)$.
\end{lemma}
\begin{proof}
We may assume that $\gamma(s)$ is an affine arc-length parameterization of the curve $\gamma$, $P_{i}=\gamma(0)=(0,0)$, $\gamma'(0)=(1,0)$ and $\gamma''(0)=(0,1)$. Then $\gamma'''(0)=(-\mu(0),0)$ and thus
we write $\gamma(s)=(s,s^2/2)+(o(s^2),o(s^3))$ with $P_{i-1}=\gamma(-s)$, $P_{i+1}=\gamma(s+o(s^2))$. Denoting $\gamma(t)=P_{i+2}$, we must show that $t=2s+o(s^2)$.
Write
\begin{eqnarray*}
P_{i+2}-P_{i-1}&=&\left(t+s+o(s^2),\frac{t^2-s^2}{2}+o(s^3)\right)\\
&=&(t+s)\left( 1+o(s),\frac{t-s}{2}+o(s^2)  \right).
\end{eqnarray*}
Also, 
\begin{eqnarray*}
P_{i+1}-P_{i}&=&\left( s+o(s^2),\frac{s^2}{2}+o(s^3)\right)\\
&=&s\left(1+o(s),\frac{s}{2}+o(s^2)\right).
\end{eqnarray*}
Since $P_{i+2}-P_{i-1}\parallel P_{i+1}-P_i$, we have
\begin{equation}\label{eq:tequal2s}
(1+o(s))(\frac{s}{2}+o(s^2))=(1+o(s))(\frac{t-s}{2}+o(s^2)).
\end{equation}
Thus $t=2s+o(s)$. Using this result in \eqref{eq:tequal2s} we obtain that in fact \linebreak $t=2s+o(s^2)$, thus proving the lemma.
\end{proof}

\begin{corollary}\label{Cor:UniformlySampling}
For a closed convex curve $\gamma$ with affine length $L$, let $\{P_1,...,P_m\}$ be the trapezoidal polygon obtained by the algorithm described 
in the beginning of this section, with $d(P_1,P_2)=d(P_2,P_3)=\frac{L}{n}$.
Denote by $\{P_1,...,{\overline P}_4,...,{\overline P}_n\}$ the affinely uniform sample along $\gamma$, i.e.,  $d({\overline P}_i,P_1)=(i-1)\frac{L}{n}$. Then 
$$
\lim_{n\to\infty}\frac{m}{n}=1
$$ 
and, for any $1\leq i\leq \min(m,n)$, 
 $$
\lim_{n\to\infty} d(P_i,{\overline P}_i)=0.
 $$ 
 \end{corollary}
 \begin{proof}
 The above lemma says that $d(P_{i+1},P_i)=s+o(s)$, $1\leq i\leq m$, where $s=\frac{L}{n}$. Thus
 $$
 d(P_i,P_1)=\frac{(i-1)L}{n}+o(1), 
 $$
 thus proving the corollary.
 \end{proof}
 
Corollary \ref{Cor:UniformlySampling} says that, for $n$ large, the convex equal-area polygon constructed above gives a sampling of the curve $\gamma$ that is approximately uniform with respect to affine arc-length. 

It is natural to ask if, given a closed convex curve $\gamma$, a point $P_1$ on it and $n$, there exists a convex equal-area $n$-gon inscribed in $\gamma$ with $P_1$ as a vertex. 
We believe that this is true for odd $n$, but not for even $n$. We plan to consider this question in a future work.

\section{ Discrete planar affine geometry and the six sextactic edges theorem}\label{section:DiscreteAffineGeometry}

\subsection{Affine curvature and sextactic edges}

At each vertex, one defines the affine normal vector $\n_i$ by
$$
\n_i=\v_{i+\frac{1}{2}}-\v_{i-\frac{1}{2}}=P_{i-1}+P_{i+1}-2P_{i}.
$$
If we assume that the polygon is convex equal-area, the determinants $[\v_{i+\frac{1}{2}},\n_i]=[\v_{i-\frac{1}{2}},\n_i]$ equal some constant $l$, independent of $i$. 
By applying an affine transformation, we may assume that $l=1$, and from now on we shall assume so.

Note that  $\n_{i+1}-\n_i$ is parallel to $\v_{i+\frac{1}{2}}$ and 
\begin{equation}\label{def:AffineCurvature}
\n_{i+1}-\n_i=-\mu_{i+\frac{1}{2}}\v_{i+\frac{1}{2}},
\end{equation}
where $\mu_{i+\frac{1}{2}}$ is defined by \eqref{RatioDiagonalSide}. We shall call $\mu_{i+\frac{1}{2}}$ the {\it discrete affine curvature} of the edge $\v_{i+\frac{1}{2}}$.
Let $\Delta\mu(i)=\mu_{i+\frac{1}{2}}-\mu_{i-\frac{1}{2}}$. A {\it sextactic edge} is an edge $\v_{i+\frac{1}{2}}$ such that $\Delta\mu(i)\cdot\Delta\mu(i+1)\leq 0$.

\subsection{ The six sextactic edges theorem}

In this section we prove a discrete analog of the six vertices theorem (\cite{Buchin83},\cite{Ovsienko01}). We begin with the following lemma:

\begin{lemma}\label{quadratic}
Consider a convex equal-area polygon $\{P_1,...,P_n\}$, $P_i=(x_i,y_i)$. Then, for any quadratic function $q(x,y)$, 
$$
\sum_{i=1}^n\Delta\mu(i)q(x_i,y_i)=0
$$
\end{lemma}

\begin{proof}
The case $q(x,y)=1$ is trivial. Consider now $q(x,y)=x$, the case $q(x,y)=y$ being analogous. Denote $\n_i=(n_i^1,n_i^2)$. A straightforward calculation shows that
\begin{eqnarray*}
\sum_{i=1}^n [\mu_{i+\frac{1}{2}}-\mu_{i-\frac{1}{2}}]\ x_i&=& -\sum_{i=1}^n\mu_{i+\frac{1}{2}}(x_{i+1}-x_i)\\
&=& \sum_{i=1}^n(n_{i+1}^1-n_i^1)=0.
\end{eqnarray*}

Now consider $q(x,y)=x^2$, which is similar to the case $q(x,y)=y^2$.
\begin{eqnarray*}
&&\sum_{i=1}^n [\mu_{i+\frac{1}{2}}-\mu_{i-\frac{1}{2}}]\ x_i^2=-\sum_{i=1}^n\mu_{i+\frac{1}{2}}(x_{i+1}^2-x_i^2)\\
&=& -\sum_{i=1}^n\mu_{i+\frac{1}{2}}(x_{i+1}-x_i)(x_{i+1}+x_i)\\
&=& \sum_{i=1}^n(n_{i+1}^1-n_i^1)(x_{i+1}+x_i)\\
&=& -\sum_{i=1}^n n_i^1\left((x_{i+1}+x_i)-(x_{i}+x_{i-1})\right)\\
&=& -\sum_{i=1}^n\left((x_{i+1}-x_i)-(x_{i}-x_{i-1})\right)\left((x_{i+1}-x_i)+(x_{i}-x_{i-1})\right)\\
&=&-\sum_{i=1}^n\left( (x_{i+1}-x_i)^2-(x_{i}-x_{i-1})^2    \right)=0.
\end{eqnarray*}

Since a rotation of the plane  preserves the convex equal-area property  we do not need to consider the term in $xy$, and so the lemma is proved.
\end{proof}

\begin{theorem}\label{SextaticEdges}
Any convex equal-area $n$-gon, $n\geq 6$, admits at least six sextactic edges.
\end{theorem}
\begin{proof}
Suppose by contradiction that $\Delta\mu(i)$ changes sign four times or less. 
Then there exists a quadratic function $q(x,y)$ that is positive 
in a region that contains the vertices where $\Delta\mu(i)$ is positive and negative in a region that contains the vertices where $\Delta\mu(i)$ is
negative. In fact, if there no changes of sign, just take $q=constant$. If there are just two changes of sign, take $q$ to be a linear function 
whose zero line divides the vertices with positive $\Delta\mu$ from the vertices with negative $\Delta\mu$. In the case of four changes of sign, 
consider lines $l_1$ and $l_2$ passing through edges where $\Delta\mu$ changes sign and whose intersection occurs inside the polygon. Then take 
$q$ as the product of linear functions whose zero lines are $l_1$ and $l_2$.  
The existence of such a quadratic function contradicts lemma \ref{quadratic} and so the theorem is proved.
\end{proof}

\section{Affine evolutes, parallel polygons and the affine distance symmetry set}\label{section:Evolutas}

\subsection{Parallel polygons and the affine evolute}

The affine normal vectors $\n_i$ generate affine normal lines $P_i(\lambda)=P_i+\lambda\n_i$, $\lambda\in\R$. The polygon $P(\lambda)$ whose vertices are $P_i(\lambda)$ is called the {\it $\lambda$-parallel polygon}.
The edges $\v_{i+\frac{1}{2}}(\lambda)=P_{i+1}(\lambda)-P_i(\lambda)$ of $P(\lambda)$ are parallel to the edges $\v_{i+\frac{1}{2}}$ of $P$, and in fact
\begin{equation}\label{eq:parallel}
\v_{i+\frac{1}{2}}(\lambda)=(1-\lambda\mu_{i+\frac{1}{2}})\ \v_{i+\frac{1}{2}}.
\end{equation}

\begin{figure}[htb]
 \centering
 \includegraphics[width=0.7\linewidth]{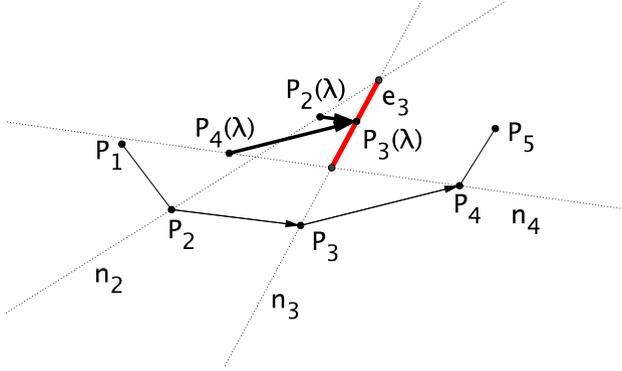}
 \caption{A cusp of a parallel at the normal line at $P_3$. The segment $\e_3$ is part of the affine evolute.}
\label{ParallelCusp}
\end{figure}

A vertex $P_i(\lambda)$ of the parallel polygon is a {\it cusp} if the determinants $[\v_{i-\frac{1}{2}},\v_{i+\frac{1}{2}}]$ and $[\v_{i-\frac{1}{2}}(\lambda),\v_{i+\frac{1}{2}}(\lambda)]$ 
have different signs. By equation \eqref{eq:parallel}, this is equivalent to 
\begin{equation}\label{eq:parallelcusp}
(1-\lambda\mu_{i-\frac{1}{2}})\cdot(1-\lambda\mu_{i+\frac{1}{2}})<0.
\end{equation}
This condition is equivalent to
$\v_{i-\frac{1}{2}}(\lambda)$ and $\v_{i+\frac{1}{2}}(\lambda)$ being at the same side relative to the normal line $P_i(\lambda),\ \lambda\in\R$ (see figure \ref{ParallelCusp}).

For each pair $(i,i+1)$, define the node $Q_{i+\frac{1}{2}}$ as the intersection of the normal lines at $P_i$ and $P_{i+1}$. Then 
$$
Q_{i+\frac{1}{2}}=P_i+\frac{1}{\mu_{i+\frac{1}{2}}}\n_i=P_{i+1}+\frac{1}{\mu_{i+\frac{1}{2}}}\n_{i+1}.
$$
Let $\e_i$ denote the subset of the normal line $P_i(\lambda)$ with boundary $\{Q_{i-\frac{1}{2}},Q_{i+\frac{1}{2}}\}$ that does not contain  $P_i$. Observe that if $\mu_{i-\frac{1}{2}}\mu_{i+\frac{1}{2}}>0$, then $\e_i$
is the segment $\overline{Q_{i-\frac{1}{2}}Q_{i+\frac{1}{2}}}$, while if $\mu_{i-\frac{1}{2}}\mu_{i+\frac{1}{2}}<0$, then $\e_i$ is the complement of this segment. If $\mu_{i-\frac{1}{2}}$ or $\mu_{i+\frac{1}{2}}$ is zero, then $\e_i$
is a half-line. The graph with nodes $Q_{i+\frac{1}{2}}$ and edges $\e_i$ is called the {\it affine evolute} of $P$. Observe that the affine evolute is a bounded polygon if and only if $\mu_{i+\frac{1}{2}}>0$, for any $1\leq i\leq n$.
From equation \eqref{eq:parallelcusp} one can easily prove the following proposition, which is well-known for smooth curves:

\begin{proposition}\label{ParallelEvolute}
Every cusp of a parallel belong to the affine evolute.
\end{proposition}

The following proposition suggests that affinely regular polygons are the discrete counterparts of ellipses:

\begin{proposition}\label{PropPolRegular}
The affine evolute of a convex equal-area polygon reduces to a point if and only if the polygon is affinely regular.
\end{proposition}
\begin{proof} 
If the polygon is affinely regular, then it is obvious that the affine evolute reduces to a point. Now suppose that the affine evolute of $P$ reduces to a point.
Since $Q_{i-\frac{1}{2}}=Q_{i+\frac{1}{2}}$, we obtain 
$$
P_i+\frac{1}{\mu_{i+\frac{1}{2}}}\n_i=P_i+\frac{1}{\mu_{i-\frac{1}{2}}}\n_i.
$$
Thus $\mu_{i+\frac{1}{2}}=\mu_{i-\frac{1}{2}}$. Since this holds for any $i$, we conclude that $\mu_{i+\frac{1}{2}}$ is independent of $i$. Now theorem 2 of \cite{Harel03} implies that $P$ is affinely regular.
\end{proof}

\subsection{ Cusps of the affine evolute}

\begin{figure}[htb]
 \centering
 \includegraphics[width=0.7\linewidth]{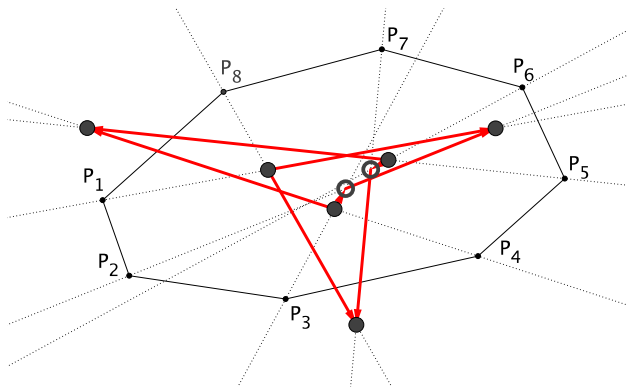}
 \caption{An oriented affine evolute. The filled dots are its cusps, while the unfilled dots are its ordinary vertices.}
\label{EvolutaAfim2}
\end{figure}

We now give an orientation to the edges of the affine evolute. 
The orientation of $\e_i$ is defined as the orientation of $\n_i$. 
A {\it cusp}  of the affine evolute is a node $Q_{i+\frac{1}{2}}$ whose adjacent pair of edges is not coherently oriented (see figure \ref{EvolutaAfim2}). 

\begin{proposition}\label{SextaticCusps}
$Q_{i+\frac{1}{2}}$ is a cusp of the affine evolute if and only if $\v_{i+\frac{1}{2}}$ is a sextactic edge of the polygon.
\end{proposition}
\begin{proof}
We shall assume $\mu_{i-\frac{1}{2}}>0$, $\mu_{i+\frac{1}{2}}>0$ and $\mu_{i+\frac{3}{2}}>0$, the other cases being analogous. Since
\begin{equation*}
Q_{i+\frac{1}{2}}-Q_{i-\frac{1}{2}}=\frac{\Delta\mu(i)}{\mu_{i-\frac{1}{2}}\mu_{i+\frac{1}{2}}}\n_i,
\end{equation*}
$\e_i$ is oriented coherently with $Q_{i+\frac{1}{2}}-Q_{i-\frac{1}{2}}$ if and only if $\Delta\mu(i)>0$. The same holds for $\e_{i+1}$ and thus 
$Q_{i+\frac{1}{2}}$ is a cusp if and only if $\Delta\mu(i)\Delta\mu(i+1)<0$, thus proving the proposition. 
\end{proof}
 
\begin{figure}[htb]
 \centering
 \includegraphics[width=0.6\linewidth]{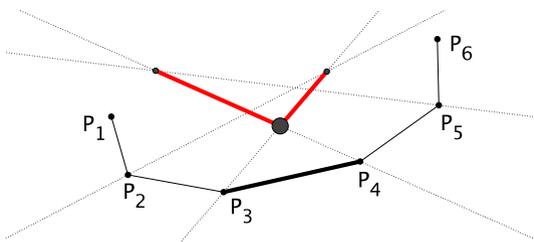}
 \caption{Correspondence between the cusp $Q_{3+\frac{1}{2}}$ of the affine evolute (big dot) and the sextactic edge $\v_{3+\frac{1}{2}}$ of the polygon.  In the figure, $\mu_{3+\frac{1}{2}}>\mu_{2+\frac{1}{2}}>0$ and $\mu_{3+\frac{1}{2}}>\mu_{4+\frac{1}{2}}>0$.}
\label{EvolutaAfim4}
\end{figure}

Following \cite{Tabach00}, the normal lines $P_i(\lambda),\ \lambda\in\R$ form an exact system, which means that the parallel polygons are closed. Then the discrete four vertex theorem 
of \cite{Tabach00} says that the affine evolute admits at least four cusps.
The following corollary, which is a direct consequence of proposition \ref{SextaticCusps} and theorem \ref{SextaticEdges}, says that, in the context of convex equal-area polygons, the affine evolute has at least six cusps. 

\begin{corollary}\label{cor:evolutecusps}
The affine evolute of a convex equal-area polygon has at least six cusps.
\end{corollary}

\subsection{  Affine distance symmetry set }

\begin{figure}[htb]
 \centering
 \includegraphics[width=0.8\linewidth]{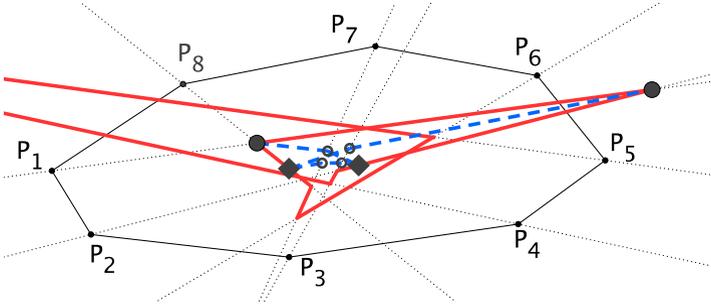}
 \caption{The dashed line is one branch of the ADSS. Observe the endpoints of the ADSS (round dots) at the cusps of the affine evolute and the cusps of the ADSS (square dots) at the affine evolute. 
 This ADSS has two other branches.}
\label{ADSS}
\end{figure}

The set of self-intersections of the parallel curves is  called {\it affine distance symmetry set} (ADSS). More precisely, the ADSS is the closure of the set of points $X$ such that 
there exist two edges 
$\v_{i+\frac{1}{2}}$ and $\v_{j+\frac{1}{2}}$, $|i-j|\geq 2$, and $\lambda\in\R$ such that the edges $\v_{i+\frac{1}{2}}(\lambda)$ and $\v_{j+\frac{1}{2}}(\lambda)$ of $P(\lambda)$ intersect 
at $X$ (see figure \ref{ADSS} and applet \cite{Applets}).

\begin{figure}[htb]
 \centering
 \includegraphics[width=0.5\linewidth]{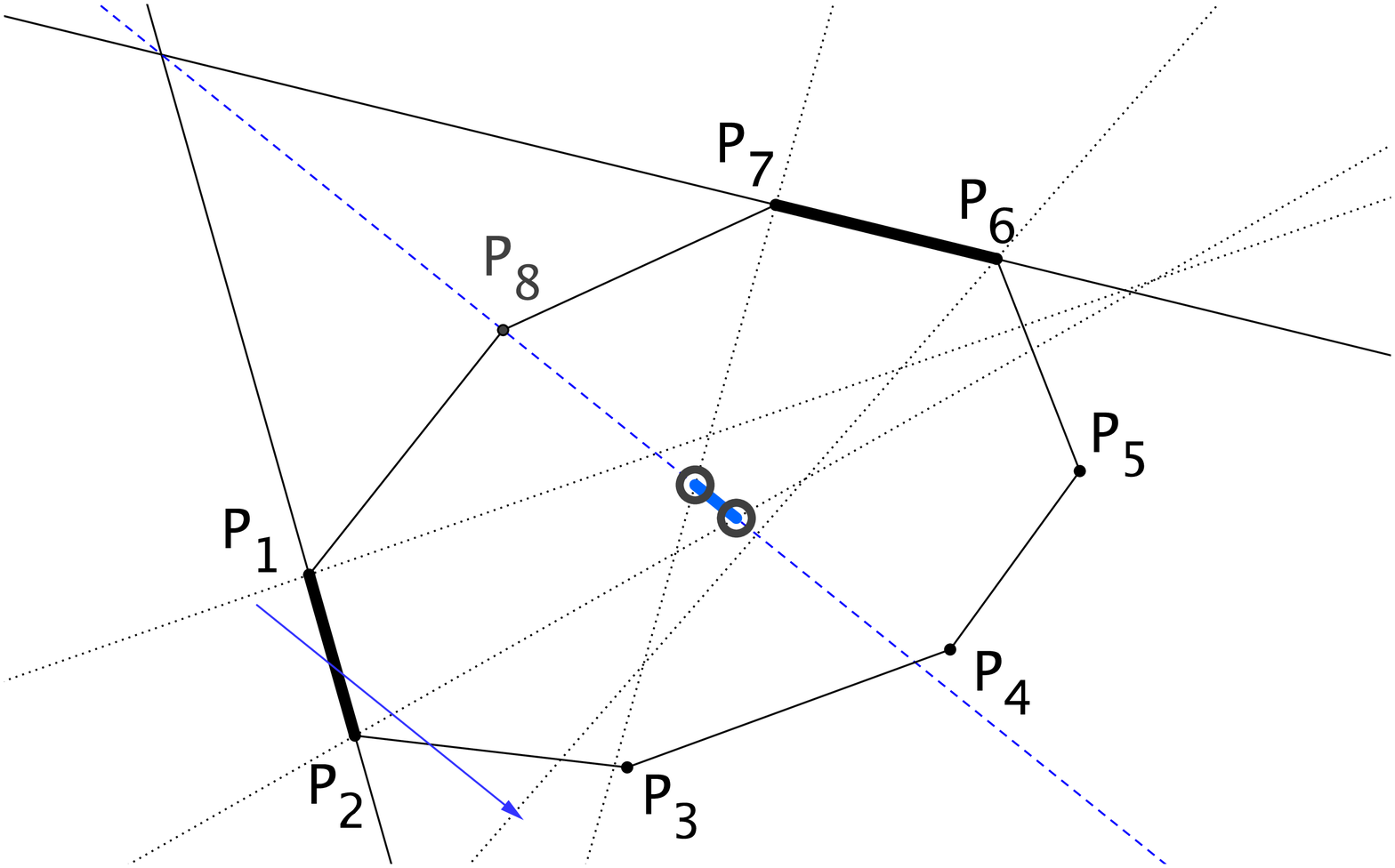}
 \caption{The dashed line is the parallel to the vector $P_6-P_7+P_2-P_1$ through the intersection of $P_1P_2$ and $P_6P_7$. This line is not the normal at $P_8$, although, by the equal-area property, $P_8$ belongs to it.
 The segment $\h(1+\frac{1}{2},6+\frac{1}{2})$ have endpoints at the normal lines at $P_7$ and $P_2$. }
\label{ADSSedge}
\end{figure}

Denote by $\h(i+\frac{1}{2},j+\frac{1}{2})$ the points of the ADSS associated with the edges  
$\v_{i+\frac{1}{2}}$ and $\v_{j+\frac{1}{2}}$. 
Since a point of $\h(i+\frac{1}{2},j+\frac{1}{2})$ is affinely equidistant to the edges $\v_{i+\frac{1}{2}}$ and $\v_{j+\frac{1}{2}}$, one can verify that  $\h(i+\frac{1}{2},j+\frac{1}{2})$ is a 
segment contained in the line through the intersection 
of the sides $\v_{i+\frac{1}{2}}$ and $\v_{j+\frac{1}{2}}$ with direction given by $\v_{i+\frac{1}{2}}-\v_{j+\frac{1}{2}}$. Moreover, the endpoints of  $\h(i+\frac{1}{2},j+\frac{1}{2})$ are points of the normal lines 
at $P_i,P_{i+1},P_j$ or $P_{j+1}$ (see figure \ref{ADSSedge}). It is possible to give a coherent orientation for $\h(i+\frac{1}{2},j+\frac{1}{2})$, but we shall not need it in this paper.  

Let $M$ be a node of the ADSS. Note that $M$ must belong to a normal line at some vertex $P_i$. More specifically, either $M$ is a common vertex of $\h(i-\frac{1}{2},j+\frac{1}{2})$ and $\h(i+\frac{1}{2},j+\frac{1}{2})$, 
in which case we denote it by $M_{i,j+\frac{1}{2}}$ or M is connected to only one edge $\h(i-\frac{1}{2},i+\frac{3}{2})$, in which case we denote it simply by $M_{i+\frac{1}{2}}$. In the latter case, we say that M is an {\it endpoint} of the ADSS (see figure  \ref{EndPoints}). The following proposition is a discrete counterpart of a well-known result of smooth curves:

\begin{figure}[htb]
\centering \fsep \subfigure[ An endpoint of the ADSS corresponding to a cusp of the affine evolute.] {
\includegraphics[width=.45
\linewidth,clip =false]{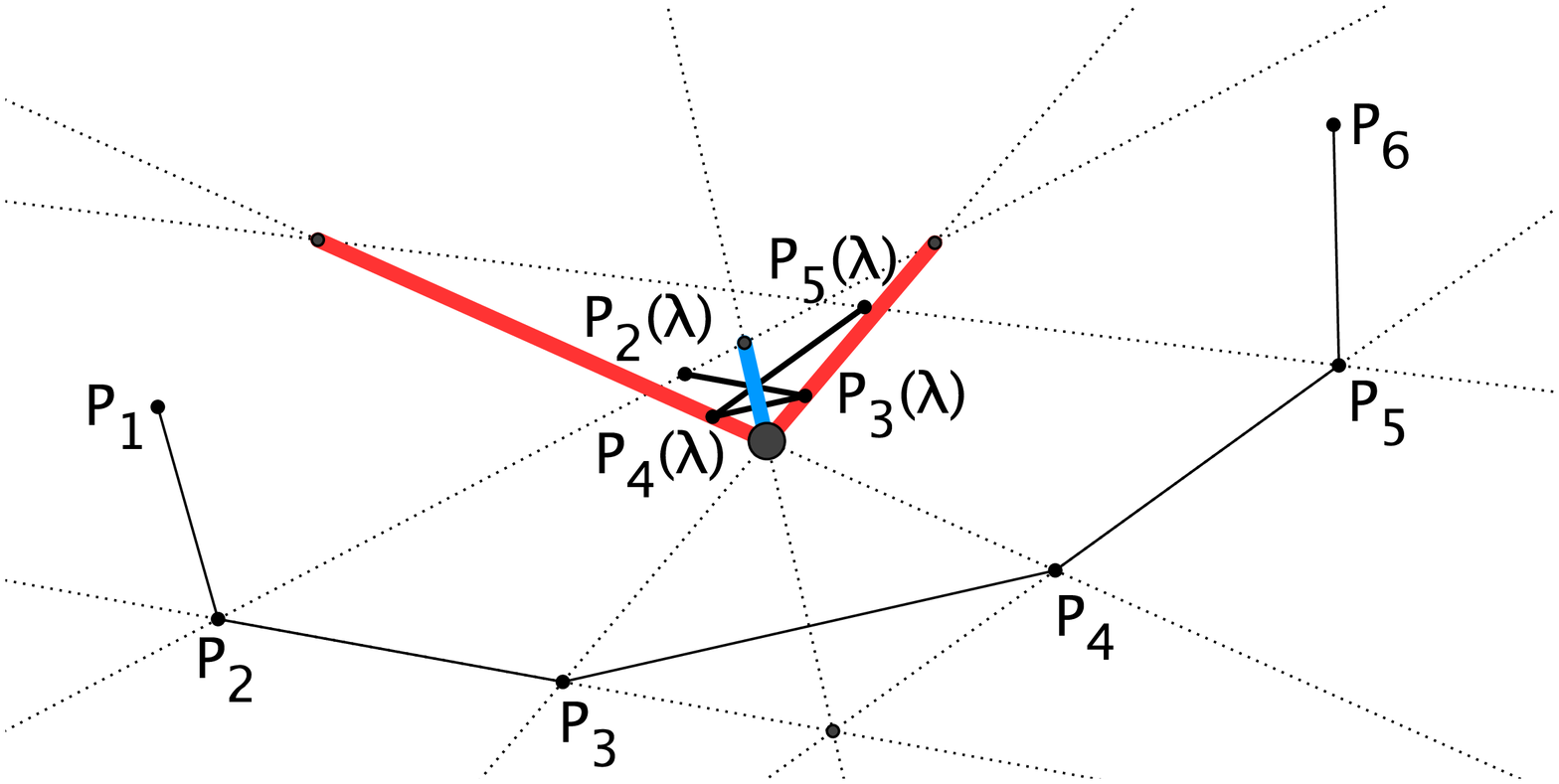}} \fsep\subfigure[
An ordinary vertex of the affine evolute that does not belong to the ADSS. ] {
\includegraphics[width=.48\linewidth,clip
=false]{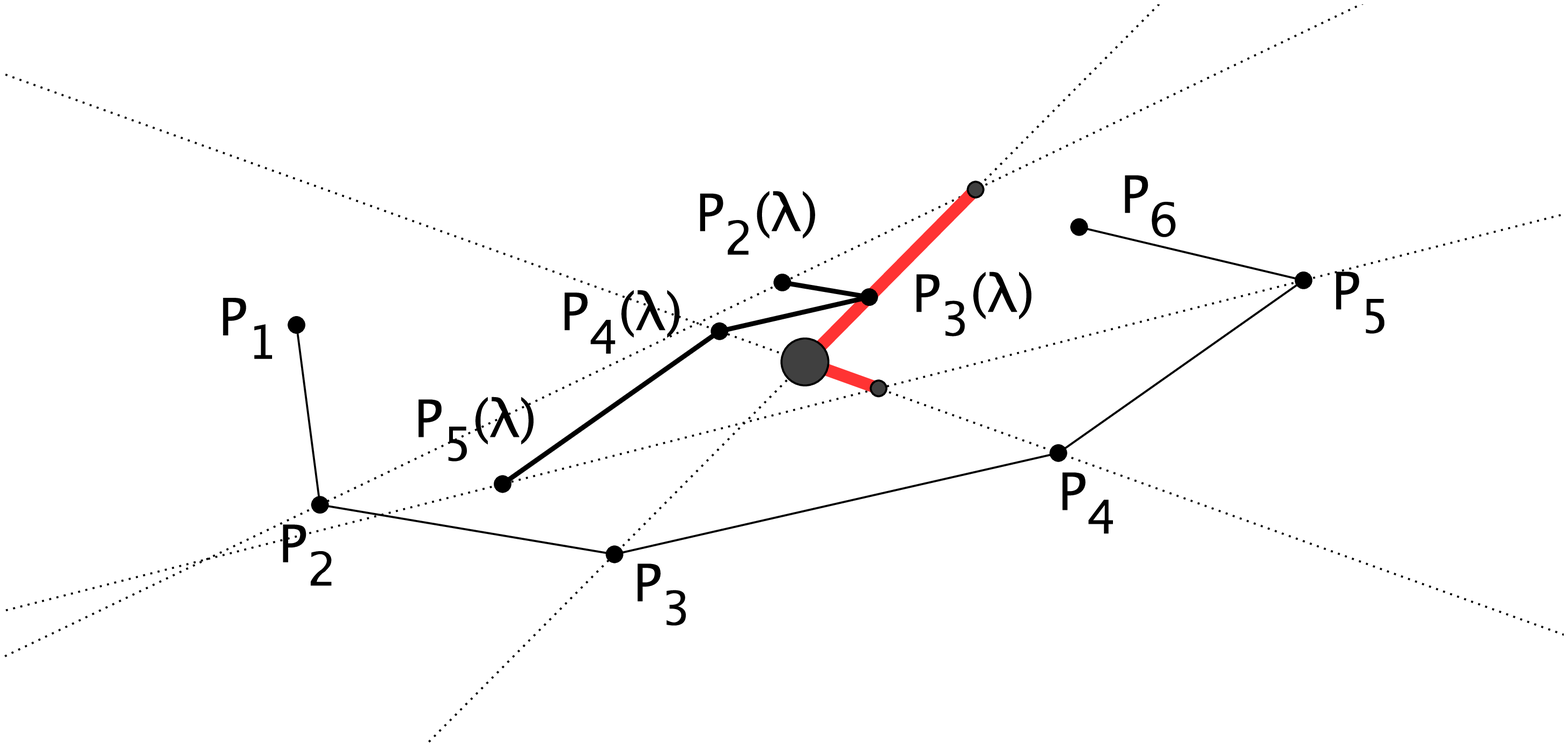}}\fsep
\caption{Endpoints of the ADSS corresponds to cusps of the affine evolute.}
\label{EndPoints}
\end{figure}

\begin{proposition}
$M_{i+\frac{1}{2}}$ is an endpoint of the ADSS if and only if $M_{i+\frac{1}{2}}$ is a cusp of the affine evolute.
\end{proposition}
\begin{proof}
$M_{i+\frac{1}{2}}$ is an endpoint of the ADSS if it is the limit of cusps of the parallels $P(\lambda)$ at the normal lines $P_i(\lambda),\ \lambda\in\R$ and $P_{i+1}(\lambda),\ \lambda\in\R$, with 
$\lambda$ converging to $(\mu_{i+\frac{1}{2}})^{-1}$.
It is now easy to see that this occurs if and only if $M_{i+\frac{1}{2}}$ is a cusp of the affine evolute. 
\end{proof}

As a consequence of the above proposition and corollary \ref{cor:evolutecusps}, we conclude that any ADSS has at least three branches. 

\bigskip

\begin{figure}[htb]
 \centering
 \includegraphics[width=0.4\linewidth]{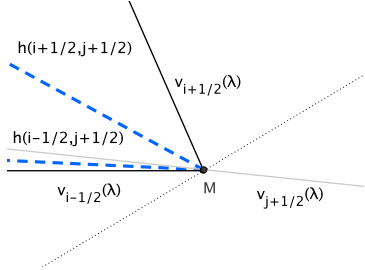}
\caption{A zoom of figure \ref{ADSS}. Observe that, in the neighborhood of a cusp $M$ of the ADSS, both $\h(i+\frac{1}{2},j+\frac{1}{2})$ and $\h(i-\frac{1}{2},j+\frac{1}{2})$ (dashed)
are on the same side of the normal at $P_i$ (dotted) as $\v_{i+\frac{1}{2}}(\lambda)$ and $\v_{i-\frac{1}{2}}(\lambda)$ (black). The remaining line in grey is $\v_{j+\frac{1}{2}}(\lambda)$, which crosses
the normal line.}
\label{ADSSZoom}
\end{figure}

Take now a node $M_{i,j+\frac{1}{2}}$ of the ADSS that is not an endpoint. 
Then $M_{i,j+\frac{1}{2}}$ is at the boundary of $\v_{i-\frac{1}{2}}(\lambda)$
and $\v_{i+\frac{1}{2}}(\lambda)$, for a certain $\lambda$, and also belongs to $\v_{j+\frac{1}{2}}(\lambda)$. We say that $M_{i,j+\frac{1}{2}}$ is a {\it cusp} of the ADSS if it is a cusp of the 
parallel $P(\lambda)$, which is equivalent to say that the edges $\v_{i-\frac{1}{2}}(\lambda)$
and $\v_{i+\frac{1}{2}}(\lambda)$ are on the same side of the normal line at $P_i$ (see figure \ref{ADSSZoom}). It follows then from proposition \ref{ParallelEvolute} that a node of the 
ADSS is a cusp if and only if belongs to the affine evolute.

\section{An affine isoperimetric inequality}\label{section:Isoperimetric}

Denote by $A=A_0$ the area bounded by the convex equal-area $n$-gon $P$. Assuming that $[\v_{i-\frac{1}{2}},\n_i]=[\v_{i+\frac{1}{2}},\n_i]=1$, the affine perimeter is defined as $L=n$. The 
following inequality is a discrete counterpart of inequality \eqref{SmoothIsoIneq}. 

\begin{theorem}\label{isoperimetric} The following isoperimetric inequality holds,
$$
\sum_{i=1}^n \mu_{i+\frac{1}{2}}\leq \frac{L^2}{2A},
$$
with equality if and only if $P$ is affinely regular.
\end{theorem}

The proof of this theorem is based on the inequality of Minkowski for mixed areas (\cite{Flanders68,Schneider93}).
Define the mixed area of two parallel $n$-gons $P$ and $P'$ by
$$
A(P,P')=\frac{1}{2}\sum_{i=1}^n[P'_i,P_{i+1}-P_i].
$$
We remark that $A(P,P')=A(P',P)$, since
\begin{eqnarray*}
\sum_{i=1}^n[P_i,P'_{i+1}-P'_i]&=&\sum_{i=1}^n[P'_{i+1},P_{i+1}-P_{i}]\\
&=&\sum_{i=1}^n[P'_i,P_{i+1}-P_{i}],
\end{eqnarray*}
where we have used the parallelism of $P_iP_{i+1}$ with $P'_{i}P'_{i+1}$ in the last equality. When $P=P$ we write $A(P)=A(P,P)$. 

Since, for small $\lambda$, $P+\lambda\n$ is convex, we can apply the Minkowski inequality for the convex parallel polygons $P$ and $P+\lambda\n$, which says that
$$
A(P,P+\lambda\n)^2\geq A(P)A(P+\lambda\n),
$$
with equality only in case $P$ and $P+\lambda\n$ are homothetic (see \cite{Schneider93}, p.321, note 1). Now 
\begin{equation*}
A(P+\lambda\n)=A(P)+2\lambda A(P,\n)+\lambda^2 A(\n)
\end{equation*}
and
\begin{equation*}
A(P,P+\lambda\n)=A(P)+\lambda A(P,\n).
\end{equation*}
Thus we conclude that
\begin{equation*}
A(P,\n)^2\geq A(P)A(\n),
\end{equation*}
which is the Minkowski inequality for $P$ and the possibly non-convex polygon $\n$. 

We can now complete the proof of theorem \ref{isoperimetric}. 
We have that
$$
A(\n)=\frac{1}{2}\sum [\n_i,\n_{i+1}-\n_i]=\frac{1}{2}\sum \mu_{i+\frac{1}{2}}[\v_{i+\frac{1}{2}},\n_i]=\frac{1}{2}\sum \mu_{i+\frac{1}{2}}.
$$
On the other hand 
$$
A(P,\n)=\frac{1}{2}\sum_{i=1}^n[\n_i,P_{i+1}-P_i]=-\frac{n}{2}=-\frac{L}{2}.
$$
Since $A(P)=A$, the Minkowski inequality for $P$ and $\n$ implies that
$$
\frac{L^2}{4}\geq A\frac{1}{2}\sum \mu_{i+\frac{1}{2}},
$$
which proves the  isoperimetric inequality. Equality holds if and only if $P$ and $P+\lambda\n$ are homothetic, which is equivalent to the affine evolute of $P$ being a single point. 
By proposition \ref{PropPolRegular}, this fact  occurs if and only if the polygon $P$ is affinely regular.


\end{document}